\input amstex
\documentstyle{amsppt}
\topmatter
\title Bernstein processes, Euclidean Quantum Mechanics and Interest Rate Models
\endtitle
\author Paul Lescot 
\endauthor
\affil 
Laboratoire de Math\'ematiques Rapha\"el Salem \\
UMR 6085 CNRS-Universit\'e de Rouen \\
Bureau M2.24 \\
Avenue de l'Universit\'e, B.P.12 \\
Technop\^ole du Madrillet \\
F76801 Saint-Etienne-du-Rouvray \\
T\'el. (+ 33) (0)2 32 95 52 24 \\
Fax (+ 33) (0)2 32 95 52 86 \\
Paul.Lescot\@univ-rouen.fr \\
http://www.univ-rouen.fr/LMRS/Persopage/Lescot \\
\endaffil
\date November 28th, 2009
\enddate
\endtopmatter
\document
\NoBlackBoxes
\TagsOnRight

\newpage
\heading{Abstract}
\endheading

We give an exposition, following joint works with J.-C. Zambrini, of the link between
euclidean quantum mechanics, Bernstein processes and isovectors for the heat equation.
A new application to mathematical finance is then discussed.

\newpage

\heading{1.Euclidean Quantum Mechanics}
\endheading

Schr\"odinger's equation for a (possibly time--dependent) potential  $V(t,q)$ :

$$
i\hbar\dsize\frac{\partial\psi}{\partial t}=-\dsize\frac{\hbar^{2}}{2m}\Delta\psi+V\psi\equiv H\psi  \,\, 
$$
on $L^{2}(\bold R^{d}, dq)$ can be written,
in space dimension $d=1$, and for $m=1$:
$$
i\hbar\dsize\frac{\partial \psi}{\partial t}=-\dsize\frac{\hbar^{2}}{2}
\dsize\frac{\partial^{2}\psi}{\partial q^{2}}+V\psi\,\, .
$$

We shall henceforth treat $\theta=\sqrt{\hbar}$ as a new parameter.

In Zambrini's Euclidean Quantum Mechanics (see e.g. \cite{1}), this equation splits into :

$$
\theta^{2}\dsize\frac{\partial \psi}{\partial t}=-\dsize\frac{\theta^{4}}{2}
\dsize\frac{\partial^{2}\psi}{\partial q^{2}}+V\psi   \tag"($\Cal C_{1}^{(V)}$)"
$$

and

$$
-\theta^{2}\dsize\frac{\partial \psi}{\partial t}=-\dsize\frac{\theta^{4}}{2}
\dsize\frac{\partial^{2}\psi}{\partial q^{2}}+V\psi          \,\, , \tag"($\Cal C_{2}^{(V)}$)"
$$
the probability density being given, not by $\psi\bar \psi$ as in the usual quantum mechanics, but by $\eta\eta_{*}$, 
$\eta$ and $\eta_{*}$ denoting respectively an everywhere strictly positive solution of \thetag{$\Cal C_{1}^{(V)}$} 
and an everywhere strictly positive solution of \thetag{$\Cal C_{2}^{(V)}$}. To these data 
is associated a \it  Bernstein process \rm
$z$, satisfying the stochastic differential equation 
$$
dz(t)=\theta dw(t)+\tilde{B}(t , z(t))dt \,\, \tag{($\Cal B$)}
$$
relatively to the canonical increasing filtration of the brownian $w$, and the stochastic differential equation

$$
d_{*}z(t)=\theta d_{*}w_{*}(t)+\tilde{B_{*}}(t , z(t))dt \tag{($\Cal B_{*}$)}
$$
relatively to the canonical decreasing filtration of another brownian  $w_{*}$, where
$$
\tilde{B}\equiv_{def} \theta^{2}\dsize\frac{\dsize\frac{\partial \eta}{\partial q}}{\eta}
$$
and
$$
\tilde{B_{*}}\equiv_{def} -\theta^{2}\dsize\frac{\dsize\frac{\partial \eta_{*}}{\partial q}}{\eta_{*}}\,\, .
$$

Setting $S=-\theta^{2} \ln(\psi)$, equation \thetag{$\Cal C_{1}^{(V)}$} becomes the Hamilton--Jacobi--Bellman equation :
$$
\dsize\frac{\partial S}{\partial t}=-\dsize\frac{\theta^{2}}{2}\dsize\frac{\partial^{2}S}{\partial q^{2}}
+\dsize\frac{1}{2}(\dsize\frac{\partial S}{\partial q})^{2}-V \,\, . \tag{$\Cal C_{3}^{(V)}$}
$$

Modulo the addition of the derivatives
$E=-\dsize\frac{\partial S}{\partial t}$
and $B=-\dsize\frac{\partial S}{\partial q}$ as auxiliary unknown functions, \thetag{$\Cal C_{3}^{(V)}$} is equivalent to the vanishing
of the following differential forms :

$$
\omega = dS+Edt+Bdq \,\, ,
$$

$$
\Omega = d\omega=dEdt+dBdq \,\, ,
$$
and
$$
\beta=(E+\dsize\frac{B^{2}}{2}-V)dqdt+\dsize\frac{\theta^{2}}{2}dBdt \,\, 
$$
on a $2$--dimensional submanifold of $\bold M=\bold R^{5}$ ($(t,q,S,E,B)$ being now considered as \it independent \rm variables).
Let then
$L=\dsize\frac{1}{2}B^{2}+V$ 
denote the \it formal lagrangian \rm, 
$$
\omega_{PC}= E dt + B dq=\omega - dS
$$ 
the \it Poincar\'e--Cartan form \rm, and
$I$ the ideal of $\Cal A=\wedge T^{*}(\bold M)$ generated by $\omega$, $d\omega$ and
$\beta$. By an \it isovector \rm we shall mean a vector field $N$ on $\bold M$ such that 

$\Cal L_{N}(I)\subseteq I$ ; because of the linearity of \thetag{$\Cal C_{1}^{(V)}$} the Lie algebra $\Cal G_{V}$ of these isovectors contains an infinite--dimensional
abelian ideal $\Cal J_{V}$ , that possesses a canonical supplement $\Cal H_{V}$.

In the free case ($V=0$) this canonical supplement has dimension $6$ and admits a natural basis, each element of which corresponds to a
symmetry of the underlying physical system.

Let $\Phi_{N}=-N(S)$ be the \it phase \rm associated to $N$, and $$D\equiv_{def}\dsize\frac{\partial}{\partial t}+B\dsize\frac{\partial}{\partial q
}
+\dsize\frac{\hbar}{2}\dsize\frac{\partial^{2}}{\partial q^{2}}$$
the \it formal Ito differential \rm along the Bernstein process $z$.
The following purely algebraic results are analogs of well--known theorems of classical
analytical mechanics:

\proclaim{Theorem 1.1}For each $N\in\Cal H_{V}$ one has :
\roster
\item $\Cal L_{N}(\omega_{PC})=d\Phi_{N}$ \,\, ;
\item $\Cal L_{N}(\Omega)=0$ \,\, ;
\item $\Cal L_{N}(L)+L\dsize\frac{dN^{t}}{dt}=D\Phi_{N}$ \,\, .
\endroster
\endproclaim

For a detailed proof see \cite{7}, and for complete calculations in the free case ($V=0$) see \cite{6}.

\newpage

\heading{2. Rosencrans' Theorem}
\endheading

Let $N$ be an \it isovector \rm (for $V=0$), let $\psi$ be a solution of $\Cal C_{1}^{(0)}$, and let 
$$S=-\theta^{2}\ln(\psi);$$ then $e^{\alpha N}$ maps  $(t,q,S,E,B)$ to $(t_{\alpha},q_{\alpha},S_{\alpha},E_{\alpha},B_{\alpha})$ ; setting $$e^{-\frac{S_{\alpha}}{\theta^{2}}}=\psi_{\alpha}(t_{\alpha},q_{\alpha})\,\, ,$$ it follows that $\psi_{\alpha}$ is also a solution of $(\Cal C_{1}^{(0)})$.
We shall denote
$$
e^{\alpha\hat{N}}:\psi \mapsto \psi_{\alpha}
$$
the associated one--parameter group ;
it is easily seen that, for
$$
N=N^{t}\dsize\frac{\partial}{\partial t}+N^{q}\dsize\frac{\partial}{\partial q}-\Phi_{N}\dsize\frac{\partial}{\partial S}
 + ...
 $$
then
 $$
 \hat{N}=-N^{t}\dsize\frac{\partial}{\partial t}-N^{q}\dsize\frac{\partial}{\partial q}+\dsize\frac{1}{\theta^{2}}\Phi_{N}. \,\, ,
 $$
and it follows that $N\mapsto -\hat{N}$ is a homomorphism of Lie algebras.

Let $\eta_{u}$ denote the solution de \thetag{$\Cal C_{1}^{(0)}$} with initial condition $u$ :
$$
\dsize\frac{\partial \eta_{u}}{\partial t}=-\dsize\frac{\theta^{2}}{2}\dsize\frac{\partial^{2}\eta_{u}}{\partial q^{2}}\,\, ,    
$$
and
$$
\eta_{u}(0,q)=u(q) \,\, .
$$

Let us set :
$$
\rho_{N}(\alpha ,t,q)=(e^{\alpha\hat{N}}\eta_{u})(t,q)
$$
and
$$
\psi^{N}(\alpha, q)\equiv_{def}\rho_{N}(\alpha , 0 , q)\,\, .
$$

Then
\proclaim{Theorem 2.1}(\cite{6}, pp.321--322)

$\psi^{N}$ satisfies :

$$
\dsize\frac{\partial\psi^{N}}{\partial \alpha}=-N^{t}(0,q)(-\dsize\frac{\theta^{2}}{2}\dsize\frac{\partial^{2}\psi^{N}}{\partial\alpha^{2}})
-N^{q}(0,q)\dsize\frac{\partial\psi^{N}}{\partial q}+\dsize\frac{1}{\theta^{2}}\Phi_{N}(0,q)\psi^{N}
$$
and
$$
\psi^{N}(0,q)=u(q)\,\, .
$$
\endproclaim

Whence :

\proclaim{Corollary 2.2}Let $N$ be chosen so that $N^{t}(0,q)=-1$, $N^{q}(0,q)=-\dsize\frac{1}{\theta^{2}}(aq+b)$ and
$\Phi_{N}(0,q)=cq^{2}+dq+f$, where $a,b,c,d,f$ denote real constants, then $$\eta_{u}^{V}(t,q)\equiv_{def}\psi^{N}(t,q)$$ satisfies the \lq\lq backwards heat equation with drift term $D(q)=aq+b$ and quadratic potential $V(q)=cq^{2}+dq+f+\dsize\frac{1}{2\hbar^{2}}(D(q))^{2}
-\dsize\frac{a}{2}$\rq\rq, corresponding to a vector potential $$A=\dsize\frac{aq+ b}{\theta^{2}} \,\, :$$ 
$$
\theta^{2} \dsize\frac{\partial \eta_{u}^{V}}{\partial t}=-\dsize\frac{\theta^{4}}{2}\dsize\frac{\partial^{2}\eta_{u}^{V}}{\partial q^{2}}
+(aq+b)\dsize\frac{\partial \eta_{u}^{V}}{\partial q}+(cq^{2}+dq+f)\eta_{u}^{V} \,\,\tag{$\Cal C_{4}^{(V)}$}
$$
and
$$
\eta_{u}^{V}(0,q)=u(q) \,\, .
$$
(In the case $D(q)=0$, the potential is given by $V(q)=cq^{2}+dq+f$ and $\eta_{u}^{(V)}$ satisfies  $\Cal C_{1}^{(V)}$ ; for the general case, cf. \cite{1}, pp.71--72).
\endproclaim

\newpage

\heading{3.The case of a linear potential}
\endheading

Here $V(q)=\lambda q$ ; it appears that :

$$
\eta_{u}^{V}(t,q)=e^{-\dsize\frac{\lambda^{2}}{6\theta^{2}}t^{3}}e^{\dsize\frac{\lambda tq}{\theta^{2}}} \eta_{u}(t,q-\lambda
\frac{t^{2}}{2}) \,\, .
$$

Then $\eta_{u}^{V}$ satisfies \thetag{$\Cal C_{1}^{(V)}$} ; the drift term can be written :
$$
\align
\tilde{B_{V}}(t,q)
&=\theta^{2} \dsize\frac{\partial}{\partial q}(\ln(\eta_{u}^{V})(t,q)) \\
&=\theta^{2}\dsize\frac{\partial}{\partial q}(-\dsize\frac{\lambda^2}{6\theta^{2}}t^{3}+\dsize\frac{\lambda tq}{\theta^{2}}+\ln(\eta_{u})(t,q-\lambda \frac{t^{2}}{2})) \\
&=\lambda t+\theta^{2}\dsize\frac{\partial}{\partial q}(\ln(\eta_{u}))(t,q-\lambda \dsize\frac{t^{2}}{2})
 \\
&=\lambda t+\tilde{B}(t,q-\lambda \dsize\frac{t^{2}}{2}) \,\, . \\
\endalign
$$

Therefore, we have :
$$
dz_{V}(t)=\theta dw(t)+\lambda t dt +\tilde{B}(t,z_{V}(t)-\lambda \dsize\frac{t^{2}}{2})dt\,\, .
$$

Let us set $y(t)\equiv_{def}z_{V}(t)-\lambda \dsize\frac{t^{2}}{2}$ ; then
$$
\align
dy(t)
&=dz_{V}(t)-\lambda t dt \\
&=\theta dw(t)+\tilde{B}(t,y(t))dt \,\, ,\\
\endalign
$$
\it i.e. \rm $y(t)$ is a Bernstein process $z(t)$ associated to solution $\eta_{u}$ of the free equation $(\Cal C_{1}^{(0)})$, and
$$
z_{V}(t)=z(t)+\lambda \dsize\frac{t^{2}}{2}\,\, .
$$

In other terms, the \lq\lq perturbation\rq\rq by a constant force $\lambda$ produces a deterministic translation by $\lambda \dsize\frac{t^{2}}{2}$, which is logical on physical grounds.

Details are given in \cite{7}.

\newpage

\heading{4.The case of a quadratic potential}
\endheading

For $V(t,q)=\dsize\frac{\omega^{2}q^{2}}{2}$, one finds :
$$
\eta_{u}^{V}(t,q)=\cosh(\omega t)^{-\frac{1}{2}}e^{\dsize\frac{\omega q^{2}}{2\theta^{2}}\tanh(\omega t)}\eta_{u}(\dsize\frac{\tanh(\omega t)}{\omega},\dsize\frac{q}{\cosh(\omega t)})\,\, .
$$

Whence
$$
\tilde{B_{V}}(t,q)=\omega q \tanh(\omega t)+\dsize\frac{1}{\cosh (\omega t)}\tilde{B}(\dsize\frac{\tanh(\omega t)}{\omega},\dsize\frac{q}{\cosh(\omega t)})\,\, .
$$

Details are exposed in \cite{9}, \S 5, and a more general formula is proved
in \cite{5}.
\newpage

\heading{5.An example with $D\neq 0$}
\endheading

Here, we take $a=\theta^{2}\beta$, and $b=c=d=f=0$.
With the notations of \cite{1}, pp. 71--72 (but, of course, replacing $\bold R^{3}$ with $\bold R$), $A=\beta q$, and $$V(q)=\dsize\frac{\beta^{2}q^{2}}{2}-\dsize\frac{\beta\theta^{2}}{2}\,\, .$$

Then $\eta_{u}^{V}$ satisfies
$$
\dsize\frac{\partial \eta_{u}^{V}}{\partial t}=-\dsize\frac{\theta^{2}}{2}\dsize\frac{\partial^{2}\eta_{u}^{V}}{\partial q^{2}}+{\beta q}\dsize\frac{\partial \eta_{u}^{V}}{\partial q}\,\, .
$$

It is easy to see that :
$$
\eta_{u}^{V}(t,q)=\eta_{u}(\dsize\frac{1}{2\beta}(e^{2\beta t}-1),e^{\beta t}q)\,\, .
$$

The drift term (cf.\cite{1}, p.72) is given by :
$$
\align
\tilde{B_{V}}(t,q)
&=\theta^{2}\dsize\frac{\partial}{\partial q}(\ln(\eta_{u}^{V})(t,q))-A(t,q) \\
&=e^{\beta t}\tilde{B}(\dsize\frac{1}{2\beta}(e^{2\beta t}-1),e^{\beta t}q) -
\beta q \,\, .\\
\endalign
$$

In particular, for $\eta=1$, one finds $\eta_{u}=1$, $\tilde{B_{V}}(t,q)=-\beta q$,
$z(t)=\theta w(t)$ and
$$
dz_{V}(t)=\theta dw(t)-\beta z_{V}(t)dt\,\, ,
$$

\it i.e. \rm $z_{V}(t)$ is an Ornstein--Uhlenbeck process, as expected.

\newpage

\heading{6.One-factor affine interest rate models}
\endheading

Such a model is characterized by the instantaneous rate $r(t)$, satisfying
$$
dr(t)=\sqrt{\alpha r(t) + \beta} \,\, dw(t) +(\phi - \lambda r(t)) \,\, dt
$$

(cf. \cite{3}).

Let us set $\tilde{\phi}=\phi + \dsize\frac{\lambda \beta}{\alpha}$ ; then :

\proclaim{Theorem 6.1} Let
$$
z(t)=\sqrt{\alpha r(t) + \beta} \,\, ;
$$

then $z(t)$ is a  Bernstein process for 

$$
\theta=\dsize\frac{\alpha}{2}\,\,
$$

and the potential

$$
V(t,q)=\dsize\frac{A}{q^{2}}+Bq^{2}
$$

where :

$$
A=\dsize\frac{\alpha^{2}}{8}(\tilde{\phi}-\dsize\frac{\alpha}{4})(\tilde{\phi}-\dsize\frac{3\alpha}{4})
$$

and

$$
B=\dsize\frac{\lambda^{2}}{8}\,\, .
$$
\endproclaim

\proclaim{Corollary 6.2} The isovector algebra $\Cal H_{V}$ associated with $V$
has dimension $6$ if and only if $A=0$ ; in the opposite case,
it has dimension $4$.
\endproclaim

The condition $A=0$ is equivalent to $\tilde{\phi}\in \{ \dsize\frac{\alpha}{4}, \dsize\frac{3\alpha}{4}\}$.
I am now able to explain this in terms of Bessel processes (see \cite{4} and \cite{5}).

\newpage

\heading{Acknowledgements}
\endheading

This work is an updated version of the text of my lecture given at ISMANS on October 24th, 2008.
I am happy to thank the organizers Alain LE MEHAUTE
and Alexandre WANG for their kind invitation and their hospitality.

\copyright{Paul Lescot, November 2009} 

\enddocument
\bye